\documentclass[10pt,a4paper]{amsart}
\usepackage{latexsym}
\usepackage{amsgen,amsmath,amsxtra,amssymb,amsfonts,amsthm}
\usepackage[latin1]{inputenc}

\newcommand{\Om}{\Omega}
\newcommand{\R}{\mathbb{R}}

\theoremstyle{plain}
\newtheorem{tma}{Theorem}
\newtheorem{corolario}[tma]{Corollary}
\newtheorem{lemma}[tma]{Lemma}

\theoremstyle{definition}
\newtheorem{defi}[tma]{Definition}

\theoremstyle{remark}
\newtheorem{remark}[tma]{Remark}

\begin{document}

\title[A mixed problem for
the infinity laplacian via Tug-of-War games]{A mixed problem for
the infinity laplacian via Tug-of-War games}

\author{Fernando Charro}
\address{Departamento de Matem{\'a}ticas, U. Autonoma de
Madrid\newline
\indent
28049 Madrid, Spain.}

\email{fernando.charro@uam.es,\ jesus.azorero@uam.es}

\author{Jesus Garc\'{\i}a Azorero}

\author{Julio D. Rossi}

\address{Departamento de Matem{\'a}tica, FCEyN, U. de Buenos Aires\newline
\indent Ciudad Universitaria, Pab I, (1428), Buenos Aires, Argentina.}

\email{jrossi@dm.uba.ar}

\thanks{Partially supported by project MTM2004-02223, MEC, Spain,
project BSCH-CEAL-UAM and project CCG06-UAM$\setminus$ESP-0340, CAM, Spain.
FC also supported by a FPU grant of MEC, Spain.
JDR partially supported by UBA X066 and CONICET, Argentina.}



\begin{abstract}
In this paper we prove that a function $ u\in\mathcal{C}(\overline{\Omega})$
is the continuous value of the Tug-of-War
game described in \cite{PSSW}
if and only if it is  the unique viscosity solution to the infinity laplacian with mixed
boundary conditions
$$
\left\{
\begin{array}{ll}
\displaystyle -\Delta_{\infty}u(x)=0\quad & \text{in}\ \Omega,\\
\displaystyle \ \frac{\partial u}{\partial n}(x)=0\quad & \text{on}\ \Gamma_N,\\
\displaystyle \ u(x)=F(x)\quad & \text{on}\ \Gamma_D.
\end{array}
\right.
$$

By using the results in \cite{PSSW}, it follows that this viscous
PDE problem has a unique solution, which is the unique {\it absolutely
minimizing Lipschitz extension} to the whole $\overline{\Omega}$ (in the
sense of \cite{Aronsson} and \cite{PSSW}) of the Lipschitz boundary data $ F:\Gamma_D\to\R $.
\end{abstract}

\maketitle

\section{Introduction}

A Tug-of-War is a two-person, zero-sum game, that is, two players
are in contest and the total earnings of one are the losses of the
other. Hence, one of them, say Player I, plays trying to maximize
his expected outcome, while the other, say Player II is trying to
minimize Player I's outcome (or, since the game is zero-sum, to
maximize his own outcome). Recently, these type of games have been
used in connection with some PDE problems, see
\cite{BEJ}, \cite{KS}, \cite{PSSW2}, \cite{PSSW}.

For the reader's convenience, let us first describe briefly the game
introduced in \cite{PSSW} by Y. Peres, O. Schramm, S. Sheffield and
D. Wilson. Consider a bounded domain $\Omega\subset\R^n$, and take
$\Gamma_D\subset\partial\Omega$ and
$\Gamma_N\equiv\partial\Omega\setminus\Gamma_D$. Let
$F:\Gamma_D\rightarrow\R$ be a Lipschitz continuous function. At an
initial time, a token is placed at a point
$x_0\in\overline\Omega\setminus\Gamma_D$. Then, a (fair) coin is
tossed and the winner of the toss is allowed to move the game
position to any $x_1\in\overline{B_\epsilon(x_0)}\cap
\overline{\Omega}$. At each turn, the coin is tossed again, and the
winner chooses a new game state  $x_k\in
\overline{B_\epsilon(x_{k-1})}\cap \overline{\Omega}$. Once the
token has reached some $x_\tau\in\Gamma_D$, the game ends and Player
I earns $F(x_\tau)$ (while Player II earns $-F(x_\tau)$). This is
the reason why we will refer to $F$ as the \emph{final payoff
function}. In more general models, it is considered  also a {\it
running payoff} $ f (x) $ defined in $ \Omega$, which represents the
reward (respectively, the cost) at each intermediate state $x$, and
gives rise to nonhomogeneous problems. We will assume throughout the
paper that $ f \equiv 0 $. This procedure yields a sequence of game
states $x_0,x_1,x_2,\ldots, x_\tau$, where every $x_k$ except $x_0$
are random variables, depending on the coin tosses and the
strategies adopted by the players.

Now we want to give a precise definition of the {\it value of the game}.
To this end we have to introduce some notation and put the game into its normal or strategic form
(see \cite{PSSW2} and \cite{sorin}).
The initial state $x_0\in\overline\Omega\setminus\Gamma_D$ is known to both players
(public knowledge). Each player $i$ chooses an \emph{action} $a_{0}^i\in \overline{B_\epsilon(0)}$ which is announced to the other player;
this defines an action profile $a_0=\{a_{0}^1,a_{0}^2\}\in \overline{B_\epsilon(0)}
\times \overline{B_\epsilon(0)}$. Then, the
new state $x_1\in \overline{B_\epsilon(x_0)}$ (namely, the current state plus the action)
is selected according to a probability  distribution $p(\cdot|x_0,a_0)$ in $\overline\Omega$ which, in our case, is given by the fair coin toss.
At stage $k$, knowing the history $h_k=(x_0,a_0,x_1,a_1,\ldots,a_{k-1},x_k)$,
(the sequence of states and actions up to that stage), each player $i$ chooses an
action $a_{k}^i$. If the game terminated at time $j < k$, we set $x_m = x_j$ and
$a_m=0$ for $j\leq m\leq k$.  The current state $x_k$ and the profile $a_k=\{a_{k}^1,a_{k}^2\}$
determine the distribution $p(\cdot|x_k,a_k)$ (again given by the fair coin toss) of the new state $x_{k+1}$.

Denote $H_k=(\overline\Omega\setminus\Gamma_{D})\times \big( \overline{B_\epsilon(0)}\times
\overline{B_\epsilon(0)}\times\overline\Omega\big)^k$, the set of
\emph{histories up to stage $k$}, and by $H=\bigcup_{k\geq1}H_k$ the set of all histories.
Notice that $H_k$, as a product space, has a measurable structure.
The \emph{complete history space} $H_\infty$ is the set of plays
defined as infinite sequences $(x_0,a_0,\ldots,a_{k-1},x_k,\ldots)$
endowed with the product topology. Then, the final payoff for Player I, i.e. $F$,
induces a Borel-measurable function on $H_\infty$.
A \emph{pure strategy} ${S}_i=\{S_{i}^k\}_k$ for Player $i$,  is a sequence of mappings from histories
to actions, namely, a mapping from $H$ to $\overline{B_\epsilon(0)}$ such that
${S}_{i}^{k}$ is a Borel-measurable mapping from $H_k$ to $\overline{B_\epsilon(0)}$
that maps histories ending with $x_k$ to elements of $\overline{B_\epsilon(0)}$
(roughly speaking, at every stage the strategy gives the next movement for the player,
provided he win the coin toss, as a function of the current state and the past history).
The initial state $x_0$ and a profile of strategies $\{S_I,S_{II}\}$
define (by Kolmogorov's extension theorem) a unique probability $\mathbb{P}_{S_I,S_{II}}^{x_0}$
on the space of plays $H_\infty$. We denote by $\mathbb{E}_{S_I,S_{II}}^{x_0}$ the corresponding expectation.

Then, if $S_I$ and $S_{II}$ denote the strategies adopted by
Player I and II respectively, we define the expected payoff for player I as
\[
V_{x_0,I}(S_I,S_{II})=
\left\{
\begin{split}
&\mathbb{E}_{S_I,S_{II}}^{x_0}[F(x_\tau)],\quad\text{if the game terminates a.s.}\\
&-\infty,\quad\text{otherwise.}
\end{split}
\right.
\]
Analogously, we define the expected payoff for player II as
\[
V_{x_0,II}(S_I,S_{II})=
\left\{
\begin{split}
&\mathbb{E}_{S_I,S_{II}}^{x_0}[F(x_\tau)],\quad\text{if the game terminates a.s.}\\
&+\infty,\quad\text{otherwise.}
\end{split}
\right.
\]
Finally, we can define the \emph{$\epsilon$-value of
the game for Player I} as
\[
u_{I}^\epsilon(x_0)=\sup_{S_I}\inf_{S_{II}}\,V_{x_0,I}(S_I,S_{II}),
\]
while the $\epsilon$-value of the game for Player II is defined as
\[
u_{II}^\epsilon(x_0)=\inf_{S_{II}}\sup_{S_I}\,V_{x_0,II}(S_I,S_{II}).
\]
In some sense, $u_I^\epsilon(x_0),u_{II}^\epsilon(x_0)$ are the
least possible outcomes that each player expects to get when the
$\epsilon$-game starts at $x_0$. Notice that, as in \cite{PSSW}, we penalize
severely the games that never end.

If $u_{I}^\epsilon= u_{II}^\epsilon :=u_\epsilon$, we say that
\emph{the game has a value}. In \cite{PSSW} it is shown that,
under very general hypotheses, that are fulfilled in the present
setting, the $\epsilon$-Tug-of-War game has a value.

\medskip

All these $\epsilon-$values are Lipschitz functions with respect
to the discrete distance $d^\epsilon$, see
\cite{PSSW} (but in general they are not
continuous), which converge uniformly when $
\epsilon
\to 0 $. The uniform limit as $\epsilon \to 0$ of the game values
$u_\epsilon$ is called
\emph{the continuous value} of the game that we will denote by
$u$. Indeed, see \cite{PSSW}, it turns out that $u$ is a viscosity
solution to the problem
\begin{equation}\label{Dirichlet}
\left\{
\begin{array}{ll}
\displaystyle - \Delta_{\infty}u(x)=0\quad & \text{in}\
\Omega,\\[8pt]
\displaystyle \ u(x)=F(x)\quad & \text{on}\ \Gamma_D,
\end{array}
\right.
\end{equation}
where $\Delta_\infty u = |\nabla u |^{-2} \sum_{ij} u_{x_i}
u_{x_ix_j} u_{x_j}$ is the $1-$homogeneous infinity laplacian (see
Section 2 for a discussion about the actual definition at points
where $\nabla u (x) $ vanishes). Infinity harmonic functions
(solutions to $-\Delta_\infty u =0$) appear naturally as limits of
$p-$harmonic functions (solutions to $\Delta_p u = \mbox{div}
(|\nabla u|^{p-2} \nabla u) =0$) and have applications to optimal
transport problems, image processing, etc. For limits as $p\to
\infty$ for $p-$laplacian type problems we refer to
\cite{ACJ},
\cite{BBM}, \cite{EG}, \cite{GAMPR} and references therein.

\medskip

When $\Gamma_D\equiv\partial\Omega$, it is known that problem
\eqref{Dirichlet} has a unique viscosity solution, (as proved
in \cite{Jensen}; see also \cite{BarlesBusca}, \cite{Crandall}, and in a more general framework,
\cite{PSSW}). Moreover, it is the unique AMLE (absolutely minimal
Lipschitz extension) of $F:\Gamma_D\to\R$ in the sense that $Lip_U
(u) = Lip_{\partial U \cap \Omega} (u)$ for every open set $U
\subset
\overline{\Omega}\setminus\Gamma_D$. AMLE extensions
were introduced by Aronsson in \cite{Aronsson}, see the survey
\cite{ACJ} for more references and applications of this subject.

However, when $\Gamma_D \not= \partial \Omega$ the PDE problem \eqref{Dirichlet} is
incomplete, since there is a missing boundary condition on $\Gamma_N =
\partial \Omega \setminus \Gamma_D$.
Our main concern is to find the boundary condition that completes
the problem.

Assuming that $ \Gamma_N $ is regular, in the sense that the normal vector field $ \vec{n}(x) $ is well defined and continuous for all $ x \in \Gamma_N $, we prove that it is in fact the homogeneous Neumann
boundary condition
$$
\frac{\partial u}{\partial n}(x)=0 \, , \quad x \in \Gamma_N .
$$

Let us point out that no regularity  is needed on the Dirichlet part $ \Gamma_D $, but the boundary data $ F $ has to be Lipschitz continuous.

On the other hand, instead of using the beautiful and involved proof based on game theory arguments, written in \cite{PSSW}, we give an alternative proof of the property $ - \Delta_{\infty} u = 0 $ in $ \Omega $, by using direct viscosity techniques, perhaps more natural in this context. The key point in our proof is the {\it Dynamic Programming Principle}, which, in some sense, plays the role of the mean property for harmonic functions in the infinity-harmonic case. This principle turns out to be an important
qualitative property of the approximations
of infinity-harmonic functions, and  is the
main tool to construct convergent numerical
methods in this kind of problems; see \cite{Oberman}.

We have the following result,

\begin{tma} \label{teo.main.intro}
Let $u(x)$ be the continuous value of the Tug-of-War game introduced
in \cite{PSSW}. Assume that $ \partial \Omega = \Gamma_N \cup \Gamma_D $, where $ \Gamma_N $ is of class $C^1$, and $ F $ is a Lipschitz function defined on $ \Gamma_D $.

Then,

\begin{enumerate}
\item[i)] $u(x)$ is a viscosity solution to the mixed boundary value problem
\begin{equation}\label{neumann.intro}
\left\{
\begin{array}{ll}
\displaystyle -\Delta_{\infty}u(x)=0\quad & \text{in}\
\Omega,\\[8pt]
\displaystyle \ \frac{\partial u}{\partial n}(x)=0\quad & \text{on}\
\Gamma_N,\\[8pt]
\displaystyle \ u(x)=F(x)\quad & \text{on}\ \Gamma_D.
\end{array}
\right.
\end{equation}

\item[ii)] Reciprocally, assume that
$\Omega$ verifies that for every $z \in \overline{\Omega}$ and
every $x^*\in\Gamma_N$ $z \neq x^*$ that
$$
\Big \langle \frac{x^*-z}{|x^*-z|} ; n (x^*) \Big\rangle >0.
$$
Then, if $u(x)$ is a viscosity solution to \eqref{neumann.intro},
it coincides with the unique continuous value of the game.
\end{enumerate}
\end{tma}

The hypothesis imposed on $\Omega$ in part ii) holds whenever $
\Gamma_N $ is strictly convex. The first part of the theorem comes
as a consequence of the Dynamic Programming Principle read in the
viscosity sense. To prove the second part we will use that the
continuous value of the game is determined by the fact that it
enjoys comparison with quadratic functions in the sense described
in \cite{PSSW}.

\medskip

We have found a PDE problem, \eqref{neumann.intro}, which allows
to find both the continuous value of the game and the AMLE of the
Dirichlet data $F$ (which is given only on a subset of the
boundary) to $\overline{\Omega}$. To summarize, we point out that
a complete equivalence holds, in the following sense:

\begin{tma} \label{tma.equivalence}
It holds
$$
u \mbox{ is AMLE of }F
\Leftrightarrow u \mbox{ is the value of the game.} \Leftrightarrow u \mbox{ solves
}\eqref{neumann.intro}.
$$
\end{tma}

\noindent The first equivalence was proved in \cite{PSSW} and the second one
is just Theorem
\ref{teo.main.intro}.

\medskip

Another consequence of Theorem \ref{teo.main.intro} is the
following:

\begin{corolario} \label{corolario.unicidad}
There exists a unique viscosity solution to \eqref{neumann.intro}.
\end{corolario}

The existence of a solution is a consequence of the existence of a
continuous value for the game together with part i) in the
previous theorem, while the uniqueness follows by uniqueness of
the value of the game and part ii).

Note that to obtain uniqueness we have to invoke the uniqueness of
the game value. It should be desirable to obtain a direct proof
(using only PDE methods) of existence and uniqueness for
\eqref{neumann.intro} but it is not clear how to find the appropriate
perturbations near $\Gamma_N$ to obtain uniqueness (existence
follows easily by taking the limit as $p \to \infty$ in the mixed
boundary value problem problem for the $p-$laplacian).

\begin{remark}
Corollary \ref{corolario.unicidad} allows to improve the convergence
result given in \cite{GAMPR} for solutions to the Neumann problem
for the $p-$laplacian as $p\to \infty$. The uniqueness of the limit
holds under weaker assumptions on the data (for example, $\Omega$
strictly convex).
\end{remark}

\medskip

The rest of the paper is devoted to the proof of Theorem
\ref{teo.main.intro}. In Section 2 we prove
part i) of the theorem and in Section 3 we prove part ii).

\section{The continuous value of the game is a viscosity solution to the
mixed problem}

As we have already mentioned in the Introduction, it is shown in
\cite{PSSW} that the continuous value of the game $u$ is infinity harmonic within
$\Omega$ and, in the case that $\Gamma_D=\partial\Om$, it
satisfies a Dirichlet boundary condition $u=F$ on
$\partial\Omega$.

\medskip

In this paper, we are concerned with the case in which
$\partial\Omega=\Gamma_D
\cup
\Gamma_N$ with $\Gamma_N\neq\emptyset$. Our aim in the present
section is to prove that $u$ satisfies an homogeneous
Neumann boundary condition on $\Gamma_N$, namely
\begin{equation}\label{neumann}
\left\{
\begin{array}{ll}
\displaystyle -\Delta_{\infty}u(x)=0\quad & \text{in}\
\Omega,\\[8pt]
\displaystyle \ \frac{\partial u}{\partial n}(x)=0\quad & \text{on}\
\Gamma_N,\\[8pt]
\displaystyle \ u(x)=F(x)\quad & \text{on}\ \Gamma_D,
\end{array}
\right.
\end{equation}
in the viscosity sense, where
\begin{equation}\label{def.inf.lap}
\Delta_{\infty}u(x)=\left\{\begin{array}{ll}
\displaystyle \Big\langle D^2u(x)\frac{\nabla u(x)}{|\nabla u(x)|},
\frac{\nabla u(x)}{|\nabla u(x)|}\Big\rangle, \qquad & \text{if}\ \nabla u(x)\neq
0,\\[12pt]
\displaystyle \lim_{y\to x}\frac{2\big(u(y)-u(x)\big)}{|y-x|^2},
\qquad & \text{otherwise}.
\end{array}\right.
\end{equation}

In defining $\Delta_\infty u$ we have followed \cite{PSSW}. Let us
point out that it is possible to define the infinity laplacian at points
with zero gradient in an alternative way, as in \cite{bib.Juutinen}.
However, it is easy to see that both definitions are equivalent.

\medskip

To motivate the above definition, notice that $\Delta_\infty u$ is the
second derivative of $u$ in the direction of the gradient. In fact,
if $ u $ is a $C^2 $ function and we take a direction $
{v}$, then the second derivative of $ u $ in the direction of
$ {v}$ is
$$
D^2_{{v}}u(x) = \left.\frac{d^2}{dt^2} \right|_{t= 0}
u(x+t {v}) = \sum_{i,j=1}^n \frac{\partial^2 u}{\partial x_i \partial x_j}(x) v_i v_j .
$$
If $\nabla u (x) \ne {0}$, we can take
$ {v} = \dfrac {\nabla u(x)}{| \nabla u(x) |} $,
and get $ \Delta_{\infty } u(x ) = D^2_{{v}} u(x). $

In points where $\nabla u(x)=0$, no direction is preferred, and then
expression \eqref{def.inf.lap} arises from the second-order Taylor's expansion of $u$ at the point $x$,
$$
\dfrac { 2 (u(y) - u(x))}{|y-x|^2} =  \Big\langle D^2 u(x ) \dfrac { y-x}{|y-x|} ,
\dfrac { y-x}{|y-x|} \Big\rangle + o(1) .
$$
We say that, at these points, $\Delta_\infty u(x)$
is defined if $D^2u(x) $ is the same in every direction, that is,
if the limit $ \dfrac {  (u(y) - u(x))}{|y-x|^2}$ exists as $ y
\to x $.

\medskip

Because of the singular nature of \eqref{def.inf.lap} in points where $\nabla u(x)=0$,
we have to restrict our class of test functions. We will denote
\[
S(x)=\left\{\phi\in\mathcal{C}^2\ \text{near}\ x\ \text{for which}\
\Delta_\infty\phi(x)\ \text{has been defined}\right\},
\]
this is, $\phi\in S(x)$ if $\phi\in\mathcal{C}^2$ in a neighborhood
of $x$ and either $\nabla\phi(x)\neq0$ or $\nabla\phi(x)=0$ and the limit
\[
\lim_{y\to x}\frac{2\big(\phi(y)-\phi(x)\big)}{|y-x|^2},
\]
exists.

\medskip

Now, using the above discussion of the infinity laplacian, we give
the precise definition of viscosity solution to \eqref{neumann}
following \cite{Bar}.

\begin{defi} \label{def.sol.viscosa}
Consider the boundary value problem \eqref{neumann}. Then,
\begin{enumerate}
\item A lower semi-continuous function $ u $ is a viscosity supersolution if for every
$ \phi \in S(x_0)$ such that $ u-\phi $ has a
strict minimum at the point $ x_0 \in
\overline{\Omega}$ with $u(x_0)= \phi(x_0)$ we have: If $x_0\in
\Gamma_D$,
$$
F(x_0) \leq \phi (x_0);
$$
if $x_0\in\Gamma_{N}$, the inequality
$$
\max\big\{ \langle n(x_0), \nabla \phi (x_0)\rangle , \  -\Delta_\infty\phi (x_0) \big\} \ge 0
$$ holds,
and if $x_0 \in \Omega$ then we require
$$
-\Delta_\infty\phi (x_0)\ge 0,
$$
with $\Delta_\infty\phi(x_0)$ given by \eqref{def.inf.lap}.

\item An upper semi-continuous function $u$ is a subsolution if for every $ \phi  \in
S(x_0)$ such that $ u-\phi $ has a strict maximum
at the point $ x_0 \in
\overline{\Omega}$ with $u(x_0)= \phi(x_0)$ we have: If $x_0\in
\Gamma_D$,
$$
F(x_0) \geq \phi (x_0);
$$
if $x_0\in \Gamma_N$, the  inequality
$$
\min \big\{ \langle n (x_0), \nabla \phi (x_0)\rangle , \  -\Delta_\infty\phi (x_0) \big \} \le 0
$$ holds, and if
$x_0 \in \Omega$ then we require
$$
-\Delta_\infty\phi (x_0) \le 0,
$$
with $\Delta_\infty\phi(x_0)$ given by \eqref{def.inf.lap}.

\item Finally, $u$ is a viscosity solution if it is both a super- and a
subsolution.
\end{enumerate}
\end{defi}

\begin{proof}[Proof of part {\rm i)} of Theorem \ref{teo.main.intro}]
 The starting point is the following Dynamic Programming Principle, which is
 satisfied by the value of the $\epsilon-$game (see \cite{PSSW}):
\begin{equation}\label{DPP}
2 u_\epsilon(x)=\sup_{y\in
\overline{B_\epsilon(x)}\cap\bar\Omega}u_\epsilon(y) +\inf_{y\in
\overline{B_\epsilon(x)}\cap\bar\Omega}u_\epsilon(y)\qquad\forall
x\in\bar\Omega\setminus\Gamma_D,
\end{equation}
where $B_\epsilon(x)$ denotes the open ball of radius $\epsilon$ centered at $x$.

Let us check that $u$ (a uniform limit of $u_\epsilon$) is a
viscosity supersolution to \eqref{neumann}. To this end, consider
a function $\phi\in S(x_0)$ such that $u-\phi$ has a strict local
minimum at $x_0$, this is,
\[
u(x)-\phi(x)> u(x_0)-\phi(x_0),\quad x\neq x_0.
\]
Without loss of generality, we can suppose that $\phi(x_0)=u(x_0)$. Let
us see the inequality that these test functions satisfy, as a consequence of the Dynamic Programming Principle.

Let $\eta(\epsilon)>0$ such that $\eta(\epsilon) = o
(\epsilon^2)$. By the uniform convergence of $u_\epsilon$ to $u$,
there exist a sequence $x_\epsilon\rightarrow x_0$ such that
\begin{equation}\label{viscosity}
u_\epsilon(x)-\phi(x)\geq u_\epsilon(x_\epsilon)-\phi(x_\epsilon)
-\eta(\epsilon),
\end{equation}
for every $x$ in a fixed neighborhood of $x_0$.

From \eqref{viscosity}, we deduce
\[
\sup_{y\in \overline{B_\epsilon(x_\epsilon)}\cap\bar\Omega}u_\epsilon(y)
\geq\max_{y\in \overline{B_\epsilon(x_\epsilon)}\cap\bar\Omega}\phi(y)+u_\epsilon(x_\epsilon)-\phi(x_\epsilon)
-\eta(\epsilon) \] and
\[
\inf_{y\in \overline{B_\epsilon(x_\epsilon)}\cap\bar\Omega}u_\epsilon(y)
\geq\min_{y\in \overline{B_\epsilon(x_\epsilon)}\cap\bar\Omega}\phi(y)+u_\epsilon(x_\epsilon)-\phi(x_\epsilon)
-\eta(\epsilon).
\]
Then, we have from \eqref{DPP}
\begin{equation}\label{VDPP}
2 \phi(x_\epsilon)\geq\max_{y\in \overline{B_\epsilon(x_\epsilon)}
\cap\bar\Omega}\phi(y)+\min_{y\in \overline{B_\epsilon(x_\epsilon)}\cap\bar\Omega}\phi(y) - 2 \eta(\epsilon).
\end{equation}
The above expression can be read as a {\it Dynamic Programming Principle in the viscosity sense}.

It is clear that the uniform limit of $u_\epsilon$, $u$, verifies
$$
u(x) = F(x) \qquad x\in \Gamma_D.
$$

In $\overline{\Omega} \setminus \Gamma_D$ there are two possibilities:
$x_0\in\Om$ and  $x_0\in\Gamma_N$. In the former case we have to check that
\begin{equation}\label{goal.interior}
- \Delta_{\infty}\phi(x_0)\geq 0,
\end{equation}
while in the latter, what we have to prove is
\begin{equation}\label{goal.boundary}
\max\Big\{\frac{\partial\phi}{\partial {n}}(x_0), - \Delta_{\infty}\phi(x_0)\Big\}\geq0.
\end{equation}

\medskip

{\bf CASE A.}  Our aim is to prove $-\Delta_\infty\phi(x_0)\ge 0$. Notice
 that this is a consequence of the results in \cite{PSSW},
 nevertheless the elementary arguments below provide an alternative
 proof using only direct viscosity techniques.

First, assume that $x_0\in\Omega$. If $\nabla \phi (x_0) \neq 0$ we proceed
as follows.

Since $\nabla\phi(x_0)\neq0$ we also have $\nabla\phi(x_\epsilon)\neq0$ for $\epsilon$ small enough.

In the sequel, $x_1^\epsilon,x_2^\epsilon\in\bar\Omega$ will be the points such that
\[
\phi(x_1^\epsilon)=\max_{y\in \overline{B_\epsilon(x_\epsilon)}\cap\bar\Omega}\phi(y)\qquad
\text{and}\quad\phi(x_2^\epsilon)=\min_{y\in \overline{B_\epsilon(x_\epsilon)}\cap\bar\Omega}\phi(y).
\]

We remark that $x_1^\epsilon,x_2^\epsilon\in\partial B_\epsilon(x_\epsilon)$.
Suppose to the contrary that there exists a subsequence
$x_1^{\epsilon_j}\in B_{\epsilon_j}(x_{\epsilon_j})$ of maximum points of $\phi$.
Then, $\nabla\phi(x_1^{\epsilon_j})=0$ and, since $x_1^{\epsilon_j}\to x_0$ as $\epsilon_j\to0$,
we have by continuity that $\nabla\phi(x_0)=0$, a contradiction. The argument for $x_2^{\epsilon}$ is similar.

Hence, since $\overline{B_\epsilon(x_\epsilon)}\cap\partial\Omega=\emptyset$, we have
\begin{equation}\label{maximos.minimos}
x_1^\epsilon=x_\epsilon+\epsilon\left[\frac{\nabla\phi(x_\epsilon)}{|\nabla\phi(x_\epsilon)|}+o(1)\right],
\quad\text{and}\quad
x_2^\epsilon=x_\epsilon-\epsilon\left[\frac{\nabla\phi(x_\epsilon)}{|\nabla\phi(x_\epsilon)|}+o(1)\right]
\end{equation}
as $\epsilon\rightarrow0$.
This can be deduced from the fact that, for $\epsilon$ small enough $\phi$ is approximately
the same as its tangent plane.

In fact, if we write $ x_1^{\epsilon} = x_{\epsilon} + \epsilon v^{\epsilon}$ with $|v^\epsilon|=1$,
and we fix any direction $ {w} $, then the Taylor expansion of $\phi$ gives

$$
 \phi (x_{\epsilon}) + \langle \nabla \phi (x_{\epsilon}),  \epsilon v^{\epsilon}\rangle +
 o(\epsilon)=\phi (x_1^{\epsilon} ) \geq \phi (x_{\epsilon} + \epsilon  {w})
$$
and hence
$$
 \langle \nabla \phi (x_{\epsilon}),   v^{\epsilon}\rangle + o(1) \ge \frac {\phi (x_{\epsilon}
 + \epsilon {w} )- \phi(x_{\epsilon})} \epsilon = \langle \nabla \phi (x_{\epsilon}),   {w} \rangle  + o(1)
$$
for any direction $ {w} $. This implies $$ v^{\epsilon}=
\frac{\nabla\phi(x_\epsilon)}{|\nabla\phi(x_\epsilon)|}+ o(1).$$

Now, consider the Taylor expansion of second order of $\phi$
\[
\phi(y)=\phi(x_\epsilon)+\nabla\phi(x_\epsilon)\cdot(y-x_\epsilon)
+\frac12\langle D^2\phi(x_\epsilon)(y-x_\epsilon),(y-x_\epsilon)\rangle+o(|y-x_\epsilon|^2)
\]
as $|y-x_\epsilon|\rightarrow0$. Evaluating the above expansion at
the point at which $\phi$ attains its minimum in $\overline{B_\epsilon(x_\epsilon)}$, $x_2^\epsilon$, we get
\[
\phi(x_2^\epsilon)=\phi(x_\epsilon) + \nabla\phi(x_\epsilon) (x_2^\epsilon - x_\epsilon)
+\frac12 \langle D^2\phi(x_\epsilon)(x_2^\epsilon-x_\epsilon),(x_2^\epsilon-x_\epsilon)\rangle+o(\epsilon^2),
\]
as $\epsilon\rightarrow0$.

Evaluating at its symmetric point in the ball $\overline{B_\epsilon(x_\epsilon)}$,
that is given by
\begin{equation}\label{hala.madrid}
\tilde{x}_2^\epsilon = 2 x_\epsilon - x_2^\epsilon
\end{equation}
we get
\[
\phi(\tilde{x}_2^\epsilon)=\phi(x_\epsilon) - \nabla\phi(x_\epsilon) (x_2^\epsilon - x_\epsilon)
+\frac12 \langle D^2\phi(x_\epsilon)(x_2^\epsilon-x_\epsilon),(x_2^\epsilon-x_\epsilon)\rangle+o(\epsilon^2).
\]
Adding both expressions we obtain
\[
\phi(\tilde{x}_2^\epsilon) + \phi(x_2^\epsilon)- 2 \phi(x_\epsilon) =
\langle D^2\phi(x_\epsilon)(x_2^\epsilon-x_\epsilon),(x_2^\epsilon-x_\epsilon)\rangle+o(\epsilon^2).
\]
We observe that, by our choice of $x_2^\epsilon$ as the point
where the minimum is attained,
\[
\phi(\tilde{x}_2^\epsilon) + \phi(x_2^\epsilon)- 2 \phi(x_\epsilon) \le
\max_{y\in \overline{B_\epsilon(x)}\cap\bar\Omega} \phi (y)
+ \min_{y\in \overline{B_\epsilon(x)}\cap\bar\Omega} \phi (y) - 2
\phi(x_\epsilon)  \le \eta(\epsilon).
\]
Therefore
\[
0\ge
\langle D^2\phi(x_\epsilon)(x_2^\epsilon-x_\epsilon),(x_2^\epsilon-x_\epsilon)\rangle+o(\epsilon^2).
\]
Note that from \eqref{maximos.minimos} we get
$$
\lim_{\epsilon \to 0}\frac{x_2^\epsilon-x_\epsilon}{\epsilon} = -\frac{\nabla \phi}{|\nabla \phi|} (x_0).
$$
Then we get, dividing by $\epsilon^2$ and passing to the limit,
$$
0 \leq - \Delta_{\infty}\phi(x_0).
$$

Now, if $\nabla \phi (x_0) = 0$ we can argue exactly as above and moreover,
we can suppose (considering a subsequence) that
\[
\frac{(x_2^\epsilon -x_\epsilon)}{\epsilon}\to v_2\qquad\text{as}\ \epsilon\to0,
\]
for some $v_2 \in\R^n$. Thus
\[
0 \leq -
\left\langle D^2\phi(x_0)v_2,v_2\right\rangle
= - \Delta_\infty\phi(x_0)
\]
by definition, since $\phi\in S(x_0)$.

\medskip

{\bf CASE B.} Suppose that $x_0\in\Gamma_N$. There are four sub-cases to be considered
depending on the direction of the gradient $\nabla\phi(x_0)$ and the
distance of the points $x_\epsilon$ to the boundary.

\medskip

\noindent{}CASE 1:\quad{}If either $\nabla\phi(x_0)=0$, or  $\nabla\phi(x_0)\neq0$
and $\nabla\phi(x_0)\bot  n(x_0)$, then
\begin{equation}\label{eq1}
\frac{\partial\phi}{\partial {n}}(x_0)=0\quad
\Rightarrow\quad\max\Big\{\frac{\partial\phi}{\partial {n}}(x_0), -\Delta_{\infty}\phi(x_0)\Big\}\geq0,
\end{equation}
where
\[
\Delta_\infty\phi(x_0)=\lim_{y\to x_0 }\frac{2\big(\phi(y)-\phi(x_0)\big)}{|y-x_0|^2}
\]
is well defined since $\phi\in S(x_0)$.

\medskip

\noindent{}CASE 2:
\quad$\displaystyle{\liminf_{\epsilon \to 0} \frac{\text{dist}(x_\epsilon,\partial\Omega)}{\epsilon}> 1}$,
and $\nabla\phi(x_0)\neq 0$.

\medskip

Since $\nabla\phi(x_0)\neq0$ we also have $\nabla\phi(x_\epsilon)\neq0$ for $\epsilon$ small enough.
Hence, since $B_\epsilon(x_\epsilon)\cap\partial\Omega=\emptyset$, we have, as before,
\[
x_1^\epsilon=x_\epsilon+\epsilon\left[\frac{\nabla\phi(x_\epsilon)}{|\nabla\phi(x_\epsilon)|}
+o(1)\right],\quad\text{and}\quad
x_2^\epsilon=x_\epsilon-\epsilon\left[\frac{\nabla\phi(x_\epsilon)}{|\nabla\phi(x_\epsilon)|}+o(1)\right]
\]
as $\epsilon\rightarrow0$. Notice that both $x_1^\epsilon,x_2^\epsilon\to \partial B_\epsilon(x_\epsilon)$.
This can be deduced from the fact that, for $\epsilon$ small enough $\phi$
is approximately the same as its tangent plane.

Then we can argue exactly as before (when $x_0 \in \Omega$) to obtain that
\[
0 \leq - \Delta_\infty\phi(x_0).
\]

\medskip

\noindent{}CASE 3:\quad$\displaystyle{\limsup_{\epsilon \to 0}
\frac{\text{dist}(x_\epsilon,\partial\Omega)}{\epsilon} \le 1},$
and $\nabla\phi(x_0)\neq0$ points inwards $\Omega$.

\medskip

In this case, for $\epsilon$ small enough we have that $\nabla\phi(x_\epsilon)\neq0$ points inwards as well. Thus,
\[
x_1^\epsilon=x_\epsilon+\epsilon\left[\frac{\nabla\phi(x_\epsilon)}{|\nabla\phi(x_\epsilon)|}+o(1)\right]\in\Omega,
\]
while $x_2^\epsilon\in\overline\Omega\cap\overline{B_\epsilon}(x_\epsilon)$. Indeed,
$$
\frac{|x_2^\epsilon-x_\epsilon|}{\epsilon} = \delta_\epsilon\leq1.
$$

We have the following first-order Taylor's expansions,
\[
\phi(x_1^\epsilon)=\phi(x_\epsilon)+\epsilon|\nabla\phi(x_\epsilon)|+o(\epsilon),
\]
and
\[
\phi(x_2^\epsilon)=\phi(x_\epsilon)+\nabla\phi(x_\epsilon)\cdot(x_2^\epsilon -x_\epsilon)+o(\epsilon),
\]
as $\epsilon\to0$. Adding both expressions, we arrive at
\[
\phi(x_1^\epsilon)+\phi(x_2^\epsilon)-2\phi(x_\epsilon)=\epsilon|\nabla\phi(x_\epsilon)|
+\nabla\phi(x_\epsilon)\cdot(x_2^\epsilon -x_\epsilon)+o(\epsilon).
\]
Using \eqref{VDPP} and dividing by $\epsilon>0$,
\[
0\geq|\nabla\phi(x_\epsilon)|+\nabla\phi(x_\epsilon)\cdot\frac{(x_2^\epsilon -x_\epsilon)}{\epsilon}+o(1)
\]
as $\epsilon\to0$. We can write
\[
0\geq|\nabla\phi(x_\epsilon)|\cdot(1+\delta_\epsilon\cos\theta_\epsilon)+o(1)
\]
where
\[
\theta_\epsilon=\text{angle}\left(\nabla\phi(x_\epsilon),\frac{(x_2^\epsilon -x_\epsilon)}{\epsilon}\right).
\]
Letting $\epsilon\to0$ we get
\[
0\geq|\nabla\phi(x_0)|\cdot(1+\delta_0\cos\theta_0),
\]
where $\delta_0\leq1$, and
\[
\theta_0=\lim_{\epsilon\to 0}\theta_\epsilon=\text{angle}\left(\nabla\phi(x_0),v(x_0)\right),
\]
with
$$
v (x_0) = \lim_{\epsilon \to 0} \frac{x_2^\epsilon -
x_\epsilon}{\epsilon}.
$$
Since $|\nabla\phi(x_0)|\neq0$, we find out $(1+\delta_0\cos\theta_0)\leq0$, and then
 $\theta_0=\pi$ and $\delta_0=1$. Hence
\begin{equation}\label{otrominimo}
\lim_{\epsilon \to 0} \frac{x_2^\epsilon - x_\epsilon}{\epsilon} = - \frac{ \nabla \phi}{|\nabla \phi|} (x_0),
\end{equation}
or what is equivalent,
\[
x_2^\epsilon=x_\epsilon-\epsilon\left[\frac{\nabla\phi(x_\epsilon)}{|\nabla\phi(x_\epsilon)|}+o(1)\right].
\]

Now, consider $\tilde{x}_2^\epsilon=2 x_\epsilon-x_2^\epsilon$ the symmetric
point of $x_2^\epsilon$ with respect to $x_\epsilon$.  We go back to \eqref{VDPP} and use
the Taylor expansions of second order,
\[
\phi(x_2^\epsilon)=\phi(x_\epsilon) + \nabla\phi(x_\epsilon) (x_2^\epsilon - x_\epsilon)
+\frac12 \langle D^2\phi(x_\epsilon)(x_2^\epsilon-x_\epsilon),(x_2^\epsilon-x_\epsilon)\rangle+o(\epsilon^2),
\]
and
\[
\phi(\tilde{x}_2^\epsilon)=\phi(x_\epsilon) + \nabla\phi(x_\epsilon) (\tilde{x}_2^\epsilon - x_\epsilon)
+\frac12 \langle D^2\phi(x_\epsilon)(\tilde{x}_2^\epsilon-x_\epsilon),(\tilde{x}_2^\epsilon
-x_\epsilon)\rangle+o(\epsilon^2),
\]
to get
$$
\begin{array}{rl}
\displaystyle
0\ge & \displaystyle \min_{y\in \overline{B_\epsilon(x_\epsilon)}\cap\bar\Omega}\phi(y)
+\max_{y\in \overline{B_\epsilon(x_\epsilon)}
\cap\bar\Omega}\phi(y)-2 \phi(x_\epsilon) \\[8pt]
 \geq & \phi(x_2^\epsilon)+ \phi(\tilde{x}_2^\epsilon) - 2 \phi(x_\epsilon)
 \\[8pt]
 = & \displaystyle  \nabla\phi(x_\epsilon) (x_2^\epsilon - x_\epsilon)
+\nabla\phi(x_\epsilon) (\tilde{x}_2^\epsilon - x_\epsilon)
+\frac12 \langle
D^2\phi(x_\epsilon)(x_2^\epsilon-x_\epsilon),(x_2^\epsilon-x_\epsilon)\rangle
\\[8pt]
& \displaystyle +\frac12 \langle
D^2\phi(x_\epsilon)(\tilde{x}_2^\epsilon-x_\epsilon),(\tilde{x}_2^\epsilon
-x_\epsilon)\rangle+o(\epsilon^2),
\\[8pt]
= & \displaystyle
 \langle D^2\phi(x_\epsilon)(x_2^\epsilon-x_\epsilon),(x_2^\epsilon-x_\epsilon)\rangle+o(\epsilon^2),
\end{array}
$$
by the definition of $\tilde{x}_2^\epsilon$. Then, we can divide
by $\epsilon^2$ and use \eqref{otrominimo} to obtain
$$
-\Delta_{\infty}\phi(x_0) \geq 0.
$$

\medskip

\noindent{}CASE 4:\quad$\displaystyle{\limsup_{\epsilon \to 0}
\frac{\text{dist}(x_\epsilon,\partial\Omega)}{\epsilon} \le 1},$
and $\nabla\phi(x_0)\neq0$ points outwards $\Omega$.

In this case we have
\[
\frac{\partial\phi}{\partial {n}}(x_0)=\nabla\phi(x_0)\cdot n(x_0)\geq0,
\]
since $ n(x_0)$ is the exterior normal at $x_0$ and $\nabla\phi(x_0)$ points outwards $\Omega$. Thus
\[
\max\Big\{\frac{\partial\phi}{\partial {n}}(x_0), - \Delta_{\infty}\phi(x_0)\Big\}\geq0,
\]
and we conclude that $u$ is a viscosity supersolution of \eqref{neumann}.

\medskip

It remains to check that $u$ is a viscosity subsolution of \eqref{neumann}.
This fact can be proved in an analogous way, taking some
care in the choice of the points where we perform Taylor expansions.
In fact, instead of taking \eqref{hala.madrid} we have to choose
$$
\tilde{x}_1^\epsilon = 2 x_\epsilon - x_1^\epsilon,
$$
that is, the reflection of the point where the maximum in the ball
$\overline{B_\epsilon(x_\epsilon)}$ of the test function is attained.

This ends the proof. \end{proof}

\section{A solution to the mixed problem enjoys comparison with quadratic functions.}

\rm
In this section we will assume the following hypothesis on the
domain $\Omega$.

\medskip

{\bf Hypothesis} {\it For every $z \in \overline{\Omega}$ and
every $x^*\in \Gamma_N$, $z \neq x^*$ we have
$$
 \Big \langle \frac{x^*-z}{|x^*-z|} ; n (x^*) \Big \rangle >0.
$$
}
Note that this holds, for example, if $\Omega$ is strictly
convex.

\medskip

We want to prove that a viscosity solution to
\begin{equation}\label{neumann.2}
\left\{
\begin{array}{ll}
\displaystyle - \Delta_{\infty}u(x)=0 \quad & \text{in}\
\Omega,\\[8pt]
\displaystyle \frac{\partial u}{\partial n}(x)=0\quad & \text{on}\
\Gamma_N,\\[8pt]
\displaystyle u(x)=F(x)\quad & \text{on}\ \Gamma_D,
\end{array}
\right.
\end{equation}
enjoys comparison with quadratic functions from above and below.
By the results of \cite{PSSW} this turns out to be equivalent to be
the unique continuous value for the Tug-of-War game obtained as the limit of
the $u_\epsilon$.

\medskip

Let us recall the definition of comparison with quadratic
functions given in \cite{PSSW}.

\begin{defi} \label{defi.quadra}
Let $Q(r) = a r^2 + b r +c $, with $a,b,c,r \in \R$.
Let $z\in \overline{\Omega} $, we call the function
$$
\varphi (x) = Q (|x-z|)
$$
a {\it quadratic distance function}.

We say that a quadratic distance function is $*-$increasing on
$V\subset \overline{\Omega}$ if either

1) $z\not\in V$ and for every $x \in V$, we have $Q' (|x-z|) >0$, or

2) $z \in V$ and $b=0$ and $a>0$.

\medskip

Similarly, we say that a function $\varphi$ is $*-$decreasing on $V$ if
$-\varphi$ is $*-$increasing on $V$.

\medskip

\begin{enumerate}
\item We say that $u$ enjoys comparison with quadratic functions from above if
for every $V \subset \overline{V} \subset U$ and a $*-$increasing
quadratic function $\varphi$ in $V$ with quadratic term $a \le 0$
then the inequality $\varphi \ge u$ on the relative boundary
$\partial V \cap \Omega$ implies $\varphi \ge u$ in $V$.

\item Analogously, we say that
$u$ enjoys comparison with quadratic functions from below if for
every $V \subset \overline{V} \subset U$ and a $*-$decreasing
quadratic function $\varphi$ in $V$ with quadratic term $a \ge 0$
then the inequality $\varphi \le u$ on the relative boundary
$\partial V\cap \Omega$ implies $\varphi \le u$ in $V$.
\end{enumerate}
\end{defi}

We split our arguments in two lemmas.

\begin{lemma}\label{comparison.above}
If $u$ is a solution to \eqref{neumann.2} then
$u$ enjoys comparison with quadratic functions from above.
\end{lemma}

\begin{proof}
Take $\varphi$ a $*-$increasing quadratic function in
$V\subset \overline{\Omega}\setminus \Gamma_D$ with $a\le 0$ and
such that $\varphi \ge u$ on $\partial V \cap \Omega$.
We have to show that $\varphi \ge u$ in $V$.

First, observe that we can assume that $\varphi > u $ on $\partial
V \cap \Omega$. If the conclusion is valid for that kind of
functions then just take $\varphi + k$ and then the limit as $k
\to 0$ to get the conclusion for any $\varphi$ with $\varphi \ge
u$ on $\partial V \cap \Omega$. Therefore, assume that $\varphi >
u $ on $\partial V \cap \Omega$.

Now, we argue by contradiction and assume that there is a point
$x_0 \in V$ with $u (x_0) > \varphi (x_0)$. Take
$$
\varphi_\delta (x) = \varphi (x) - \delta |x-z|^2
$$
with $\delta $ small in order to have
$\varphi_\delta > u$ on $\partial V \cap \Omega$ and
$$
\max_V \ (u - \varphi_\delta) = u(x^*) - \varphi_\delta (x^*) >0.
$$
Notice that $x^* \in V \setminus (\partial V \cap \Omega)$.

We have two possibilities:

\medskip

{\bf CASE A}  If $x^* \in V \cap \Omega$, since $u$ is a viscosity subsolution to $-\Delta_\infty u=0$
in $\Omega$, we have that,
\begin{equation}\label{con.mosca}
-\Delta_\infty \varphi_\delta (x^*) \le 0.
\end{equation}
On the other hand, since $a \le 0$ we get (just by differentiation of
the explicit expression of a quadratic function)
$$
-\Delta_\infty \varphi_\delta (x^*) = -2a + 2\delta > 0.
$$
This contradicts \eqref{con.mosca}.

\medskip

{\bf CASE B} If $x^* \in V \cap \partial \Omega$, since $u$ is a viscosity subsolution to
\eqref{neumann.2}, we have that
$$
\min \left\{ \frac{\partial \varphi_\delta}{\partial n}  (x^*), \
 -\Delta_\infty \varphi_\delta (x^*) \right\} \le 0.
$$
Therefore, we have
\begin{equation}\label{trucha}
\frac{\partial \varphi_\delta}{\partial n}  (x^*) \le 0 \qquad \mbox{ or } \qquad
-\Delta_\infty \varphi_\delta (x^*) \le 0.
\end{equation}

On the other hand, using again the fact that $a\le 0$, we obtain,
as before,
\begin{equation}\label{contra.1}
-\Delta_\infty \varphi_\delta (x^*) > 0.
\end{equation}

Now, we argue as follows, since $\varphi$ is a $*-$increasing
quadratic function on $V\subset \overline{\Omega}\setminus
\Gamma_D$ (see Definition \ref{defi.quadra}) we have that either

1) $z\not\in V$ and for every $x \in V$, we have $Q' (|x-z|) >0$, or

2) $z \in V$ and $b=0$ and $a>0$.

Note that since $a\le 0$ the second case, 2), is not possible.

Therefore, in the first case, 1), we have that
\begin{equation}\label{contra.2}
\begin{array}{rl}
\displaystyle \frac{\partial \varphi_\delta}{\partial n} (x^* ) & \displaystyle =
\frac{\partial \varphi}{\partial n} (x^* ) - \delta \frac{\partial |x-z|^2}{\partial n} (x^* )
\\[8pt]
& \displaystyle=  \Big( Q' (|x^*-z|)  - 2 \delta |x^*-z|
\Big)\Big \langle
\frac{x^*-z}{|x^*-z|} , n (x^*) \Big
\rangle >0,
\end{array}
\end{equation}
choosing $\delta$ smaller if necessary. We are using here that $Q' (|x^*-z|) >0$ and that
$$
\Big  \langle \frac{x^*-z}{|x^*-z|} , n (x^*) \Big \rangle >0.
$$
Inequalities \eqref{contra.1} and \eqref{contra.2} contradict
\eqref{trucha}. The proof is now complete.
\end{proof}

\begin{remark} We have only used that $u$ is a viscosity subsolution
to \eqref{neumann.2} to prove this lemma.
\end{remark}

In an analogous way (but using only that $u$ is a viscosity supersolution) we can prove that,

\begin{lemma}\label{comparison.below}
If $u$ is a solution to \eqref{neumann.2} then
$u$ enjoys comparison with quadratic functions from below.
\end{lemma}

\begin{proof}
It is analogous to the proof of the previous lemma.
\end{proof}

\begin{remark}
It is possible to relax the geometric assumption on $\Omega$ in Lemmas \ref{comparison.above} and \ref{comparison.below} to allow convex domains with flat pieces of boundary if we assume that, for every $z \in \overline{\Omega}$ and
every $x^*\in \Gamma_N$, $z \neq x^*$, we have
\begin{equation}\label{hyp}
 \Big \langle \frac{x^*-z}{|x^*-z|} , n (x^*) \Big \rangle \geq0,
\end{equation}
and, for all $V\subset\overline{\Omega}\setminus\Gamma_D$ there exists $p\in\R^n$, $|p|=1$, such that $\langle p,n(x^*)\rangle>0$ for all $x^*\in\Gamma_N\cup V$ for which \eqref{hyp} holds with an equality. Then, the proof of the Lemmas, can be carried out as above considering the perturbation
\[
\tilde\varphi_\delta (x) = \varphi (x) - \delta |x-z|^2 + \delta^2 \langle p, (x-z)\rangle
\]
(for $p$ as in the above hypothesis) instead of $\varphi_\delta$.
\end{remark}

From these two lemmas and the results of \cite{PSSW} we can
easily deduce that any viscosity solution to \eqref{neumann.2} has to be
the unique continuous value of the Tug-of-War game.

\begin{proof}[Proof of part {\rm ii)} of Theorem \ref{teo.main.intro}]
The previous two lemmas show that a viscosity solution to
\eqref{neumann.2} has comparison with quadratic functions from
above and below, hence, by the results of \cite{PSSW}, it is the
unique continuous value of the game.
\end{proof}

\

\noindent{\bf Acknowledgments:} The authors would like to thank
Juan J. Manfredi for many useful suggestions and conversations.

\

\end{document}